\documentclass[12pt]{article}

\usepackage[pdftex]{color, graphicx}
\usepackage{setspace}
\usepackage[utf8]{inputenc}
\usepackage{amsmath}
\usepackage[pdftex,bookmarks,colorlinks]{hyperref}
\usepackage{latexsym}
\usepackage{amsfonts}
\usepackage{amssymb}
\usepackage{amsthm}
\usepackage{paralist}
\usepackage{mathrsfs}
\usepackage[super]{nth}
\begin{document}
\newtheorem{thm}{Theorem}
\newtheorem{cor}[thm]{Corollary}
\newtheorem{lem}{Lemma}
\theoremstyle{remark}\newtheorem{rem}{Remark}
\theoremstyle{definition}\newtheorem{defn}{Definition}

\title{Weighted inequalities in ergodic theory via transference}
\author{Sakin Demir\\
Agri Ibrahim Cecen University\\ 
Faculty of Education\\
\nth{3} Floor, Office C-42\\
04100 A\u{g}r{\i}, Turkey\\
e-mail: sakin.demir@gmail.com\\
ORCID:\url{https://orcid.org/0000-0002-8019-6917}
}
\maketitle
\renewcommand{\thefootnote}{}



\renewcommand{\thefootnote}{\arabic{footnote}}
\setcounter{footnote}{0}

\begin{abstract}We first extend Calder\'on's transfer principle  to weighted weak type and strong type inequalities in various different settings including vector-valued setting. We show that our results can be used to prove a large class of inequalities in ergodic theory and ergodic $H^1$ spaces. \\
\indent
The following is one of our results: Let $X$ be a totally $\sigma$-finite measure space and $U^t$ be a one-parameter group of measure-preserving transformations of $X$. Let $w$ and $v$ be two weight functions such that $v(U^tx)$ is a locally integrable function of $t$ for almost all $x$, and  let $(S_n)$ and $(P_n)$  be sequences of sublinear and semilocal operators defined on a space of locally integrable functions on $\mathbb{R}$ commuting with translations. Define
$Sf=\sup_n|S_nf|$, $S^\sharp f=\sup_n|S_n^\sharp f|,$ and $Pf=\sup_n|P_nf|$, $P^\sharp f=\sup_n|P_n^\sharp f|$\\
where  $T^{\sharp}f(\cdot )=T(f(U^{(\cdot )}x))$, $Sf=\sup_n|S_nf|$, $S^\sharp f=\sup_n|S_n^\sharp f|,$ and $Pf=\sup_n|P_nf|$, $P^\sharp f=\sup_n|P_n^\sharp f|$. \\
\indent
Let $w$ and $v$ be two non-negative measurable functions and suppose also that $w(U^tx)$ and $v(U^tx)$ are measurable in the product of $X$ with the real line, and let $(S_n)$ and $(P_n)$  be sequences of operators as above. Define
$Sf=\sup_n|S_nf|$, $S^\sharp f=\sup_n|S_n^\sharp f|,$ and $Pf=\sup_n|P_nf|$, $P^\sharp f=\sup_n|P_n^\sharp f|$.  Suppose that there exists a constant $C_1>0$ such that
$$\int_{\mathbb{R}}|Sf(t)|^pw(U^tx)\, dt\leq C_1\int_{\mathbb{R}}|Pf(t)|^pv(U^tx)\, dt$$
and $Pf(U^{(\cdot )}x)\in L^p_{v(U^{(\cdot )}x)}(\mathbb{R})$ for almost every $x$  , for all $f\in (L^1+L^\infty )(\mathbb{R})$,  $1\leq p<\infty$. Then we have
$$\int_{X}|S^{\sharp}f(x)|^pw(x)\, d\mu (x)\leq C_1\int_{X}|P^{\sharp}f(x)|^pv(x)\, d\mu (x) $$
for all  $f\in (L^1+L^\infty )(X)$ satisfying $P^{\sharp}f\in L^p_v(X)$,  $1\leq p<\infty$.\\
\indent
Also, if there exists a constant $C_2>0$ such that for every $\lambda>0$
$$\int_{\{t:|Sf(t)|>\lambda\}}w(U^tx)\, dt\leq \frac{C_2}{\lambda^p}\int_{\mathbb{R}}|Pf(t)|^pv(U^tx)\, dt$$
and $Pf(U^{(\cdot )}x)\in L^p_{v(U^{(\cdot )}x)}(\mathbb{R})$ for almost every $x$  , for all $f\in (L^1+L^\infty )(\mathbb{R})$, then for every $\lambda>0$
$$\int_{\{x:|S^{\sharp}f(x)|>\lambda\}}w(x)\, d\mu (x)\leq \frac{C_2}{\lambda^p}\int_{X}|P^{\sharp}f(x)|^pv(x)\, d\mu (x)$$
holds for all   $f\in (L^1+L^\infty )(X)$ satisfying $P^{\sharp}f\in L^p_v(X)$, $1\leq p<\infty$.
\end{abstract}

\section{Introduction}
\indent

The transfer principle has long attracted considerable interest in analysis (see, for example, \cite{ha}, \cite{naebtag}, \cite{nhasjms}, \cite{bbmcg}, \cite{apcal}, \cite{rrcgw1}, \cite{rrcgw2}, \cite{demgct}, \cite{demect}, \cite{cfpw}, and \cite{dkosz}) because of the role it plays in proving theorems in an abstract structure. When certain conditions are satisfied it provides a bridge role between two different fields of mathematics. The transfer principle is also a subject of interest in its own right since  obtaining a transfer principle between two different fields of mathematics is a problem of finding some relation between those two fields.\\
\indent
The use of the transfer principle might be the only possible way to prove certain theorems in ergodic theory, but even when an alternative proof is available, the transfer principle always provides a much shorter proof when it is applicable to the related problem.\\
\indent
Recall the Hardy-Littlewood maximal function $Hf$ defined as
$$Hf(x)=\sup_{x\in I}\frac{1}{|I|}\int_I \left|f(x-y)\right|\,dy,$$
where the supremum is taken over all intervals $I$ in $\mathbb{R}$, $|I|$ is the measure of $I$, and $f\in L^1(\mathbb{R})$ , and the ergodic maximal function $M^\ast f $ defined as 
$$M^\ast f(x)=\sup_{x\in I}\frac{1}{|I|}\int_I\left|f(U_sx)\right|\, ds,$$
where $\{U_t:-\infty <t<\infty\}$ is a one-parameter group of measure-preserving transformations on a totally $\sigma$-finite measure space  $(X,\mathscr{B} ,\mu )$, and $f\in L^1(X)$. Understanding the connection between $Hf$ and $M^\ast f$ is probably the simplest way to understand the importance of the method of transference in analysis. It is well known by the work of A. P. Calder\'on~\cite{apcal} that the weak type $(1,1)$ inequality, and strong $L^p$ inequalities for $M^\ast f$ with $1<p<\infty$ follow from the analogous results for $Hf$. \\
\indent
A. P. Calder\'on's article appears to be the first on a transfer principle between harmonic analysis and ergodic theory in the literature. His proof relies on the use of Fubini's theorem and switching variables when passing from harmonic analysis to ergodic theory. Soon after,  the scope of Calder\'on's seminal techniques has been extended in R. R. Coifman and G. Weiss~\cite{rrcgw1}, and \cite{rrcgw2} to locally compact groups. C.~Finet and P.~Wantiezi~\cite {cfpw} extended the method of R. R. Coifman and G. Weiss~\cite{rrcgw2} to Orlicz spaces with weight concerning $L^p$ norm inequalities.\\
\indent
In S. Demir~\cite{demgct} the author proved that Calder\'on's argument can be extended to be able to replace the right-hand side with the $L^p$ norm of a different operator satisfying the conditions of Calder\'on's transfer principle in the transfer of strong $L^p$ norm inequalities (see Theorem~\ref{stp}). It was also proved by the author in S. Demir~\cite{demect} that  Calder\'on's results can be extended to the vector-valued setting (see Theorem~\ref{vvctp}). \\
\indent
H. Aimer~\cite{ha} used  the method given by A. P. Calder\'on~\cite{apcal} to show that the weighted strong $L^p$ inequalities on a locally compact abelian group can be transferred. In this work we extend Calder\'on's method to show that two-weight $L^p$ norm inequalities and weighted weak type inequalities in harmonic analysis can be transferred to ergodic theory. We also show that weighted inequalities on real Hardy spaces can be transferred to the corresponding weighted inequalities on ergodic Hardy spaces again by extending  Calder\'on's method. We also apply our results to prove some new weighted inequalities in ergodic theory and ergodic $H^1$ spaces.\\
\indent
{\bf{The following setup will be used in our discussions:}}\\
\indent
Suppose that $X$ is a measure space which is totally $\sigma$-finite and $U^t$ is a one-parameter group of measure-preserving transformations of $X$. We will also assume that for every measurable function $f$ on $X$ the function $f(U^tx)$ is measurable in the product of $X$ with the real line. $T$ will denote an operator defined on the space of locally integrable functions on the real line with the following properties: the values of $T$ are continuous functions on the real line, $T$ is sublinear and commutes with translations, and $T$ is semilocal in the sense that there exists a positive number $\epsilon$ such that the support of $Tf$ is always contained in an $\epsilon$-neighborhood of the support of $f$. \\
We will associate an operator $T^\sharp$ on functions on $X$ with such an operator $T$ as follows:\\
\indent
Given a function $f$ on $X$ let
$$F(t,x)=f(U^tx).$$
If $f$ is the sum of two functions which are bounded and integrable, respectively, then $F(t,x)$ is a locally integrable function of $t$ for almost all $x$ and therefore
$$G(t,x)=T(F(t,x))$$
is a well-defined continuous function of $t$ for almost all $x$. Since for any given value $t$, $F(t,x)$ is a measurable function of $x$ and $T$ commutes with translation, the function $G(t,x)$ has the same property. Thus $g(x)=G(0,x)$ has a meaning and we define
$$T^\sharp f=g(x).$$
\indent
Recall that an operator $T$ on $L^p(X)$ is of weak type $(p,p)$ if there exists a positive constant $C$ such that
$$\mu\{x:|Tf(x)|>\lambda \}\leq \frac{C}{\lambda^p}\|f||_p^p$$
for all $f\in L^p(X)$ and $\lambda >0$.\\
\indent
$T$ is said to be of strong type $(p,p)$ if there exists a positive constant $C$ such that
$$\|Tf||_p\leq C\|f||_p$$
for  all $f\in L^p(X)$.\\
\indent
A. P. Calder\'on~\cite{apcal} proved the following results:
 \begin{thm}\label{apcal1}Let $T_n$ be a sequence of operators as above and suppose that the operator
$Sf=\sup |T_nf|$ is of strong type $(p,p)$, $1\leq p \leq \infty$. Then the same holds for the operator  $S^\sharp f=\sup |T_n^\sharp f|$ and $\|S^{\sharp}\|\leq \|S\|$.
\end{thm}
 \begin{thm}\label{apcal2}Let $T_nf=k_n\ast f$ where $k_n$ is bounded and has bounded support. Suppose that $Sf=\sup |T_nf|$ is of weak type $(p,p)$, $1<p<\infty$, and that $\int k_n(t)\, dt$ converges and $k_n\ast \phi$ converges in $L^1$, as  $n\to \infty$, for every infinitely differentiable $\phi$ with compact support and vanishing integral. Then $T_n^{\sharp}f$ converges almost everywhere in $L^p(X)$.
\end{thm}
We will need the following result  obtained by the author in S. Demir~\cite{demect}:
 \begin{thm}\label{vvctp}Let $1\leq r\leq\infty$ and $1\leq p<\infty$. Suppose that there exists a constant $C_1>0$ such that
 $$\left(\int_{\mathbb{R}}\|Sf(t)\|_{\ell^r}^p\,dt\right)^{1/p}\leq C_1\left(\int_{\mathbb{R}}\|f(t)\|_{\ell^r}^p\,dt\right)^{1/p}$$
 where $f=(f_1,f_2,f_3,\dots )$ is a sequence of functions on $\mathbb{R}$.
 Then we have
  $$\left(\int_X\|S^\sharp f(x)\|_{\ell^r}^p\,dx\right)^{1/p}\leq C_1\left(\int_X\|f(x)\|_{\ell^r}^p\,dx\right)^{1/p}$$
  where $f=(f_1,f_2,f_3,\dots )$ is a sequence of functions on $X$.
 Suppose that there exists a constant $C_2>0$ such that for all $\lambda >0$
 $$\left|\left\{t\in\mathbb{R}:\|Sf(t)\|_{\ell^r}>\lambda\right\}\right|\leq\frac{C_2^p}{\lambda^p}\int_{\mathbb{R}}\|f(t)\|_{\ell^r}^p\,dt$$
 where $f=(f_1,f_2,f_3,\dots )$ is a sequence of functions on $\mathbb{R}$.
 Then for all $\lambda >0$ we have
 $$\left|\left\{x\in X:\|S^\sharp f(x)\|_{\ell^r}>\lambda\right\}\right|\leq\frac{C_2^p}{\lambda^p}\int_X\|f(x)\|_{\ell^r}^p\,dx$$
 where $f=(f_1,f_2,f_3,\dots )$ is a sequence of functions on $X$.
\end{thm}
\indent
The following result has been obtained by the author in S. Demir~\cite{demgct}:
\begin{thm}\label{stp}Let now $(S_n)$ and $(P_n)$  be sequences of operators as above and define
$Sf=\sup_n|S_nf|$ and $S^\sharp f=\sup_n|S_n^\sharp f|,$ and $Pf=\sup_n|P_nf|$
and $P^\sharp f=\sup_n|P_n^\sharp f|$. If there is a constant $C>0$ such that
$$\|Sf\|_p\leq C\|Pf\|_p$$
for all $f\in L^p(\mathbb{R})$ satisfying $Pf\in L^p(\mathbb{R})$  , $1\leq p\leq\infty$, then we have
$$\|S^\sharp f\|_p\leq C\|P^\sharp f\|_p$$
for all $f\in L^p(X)$ satisfying $P^\sharp f\in L^p(X)$ , $1\leq p\leq\infty$.
\end{thm}
\section{Results}
\indent

In this section we extend Calder\'on transfer principle to weighted spaces with the same setup and definitions as in the introduction.\\
\indent
Our first result of this section is the following:
 \begin{thm}\label{wifep}Let $w$ and $v$ be two non-negative measurable functions and suppose also that $w(U^tx)$ and $v(U^tx)$ are measurable in the product of $X$ with the real line, and let $T_n$ be a sequence of operators as in the introduction and 
$Sf=\sup |T_nf|$. Suppose that  there exists a constant $C_1>0$ such that
$$\int_{\mathbb{R}}|Sf(t)|^pw(U^tx)\, dt\leq C_1\int_{\mathbb{R}}|f(t)|^pv(U^tx)\, dt$$
for all $f\in (L^1+L^\infty )\cap L^p_{v(U^{(\cdot )}x)}(\mathbb{R})$ for almost every $x$, $1\leq p< \infty$. Then we have
$$\int_{X}|S^{\sharp}f(x)|^pw(x)\, d\mu (x)\leq C_1\int_{X}|f(x)|^pv(x)\, d\mu (x)$$
for all  $f\in (L^1+L^\infty )\cap L^p_v(X)$,  $1\leq p< \infty$, where $S^\sharp f=\sup |T_n^\sharp f|$.\\
\indent
Also, if there exists a constant $C_2>0$ such that for every $\lambda>0$
$$\int_{\{t:|Sf(t)|>\lambda\}}w(U^tx)\, dt\leq \frac{C_2}{\lambda^p}\int_{\mathbb{R}}|f(t)|^pv(U^tx)\, dt$$
holds for all  $f\in (L^1+L^\infty )\cap L^p_{v(U^{(\cdot )}x)}(\mathbb{R})$  for almost every $x$, $1\leq p< \infty$, then for every $\lambda>0$
$$\int_{\{x:|S^{\sharp}f(x)|>\lambda\}}w(x)\, d\mu (x)\leq \frac{C_2}{\lambda^p}\int_{X}|f(x)|^pv(x)\, d\mu (x)$$
holds for all   $f\in (L^1+L^\infty )\cap L^p_v(X)$, $1\leq p< \infty$.
\end{thm}
\begin{proof}We adapt the argument of A. P. Calder\'on~\cite{apcal} to prove our theorem.\\
 Without loss of generality we may assume that the sequence $T_n$ is finite, for if the theorem is established in this case, the general case follows by a passage to the limit. Under this assumption the operator $S$ has the same properties as the operator $T$ above.
We note that 
$$F(t,U^sx)=F(t+s,x),$$
which means that for any two given values $t_1$, $t_2$ of $t$, $F(t_1,x)$ and $F(t_2,x)$ are equimeasurable functions of $x$. On the other hand, due to translation invariance of $S$, the function $G(t,x)$ has the same property. In fact we have
$$G(t,U^sx)=S(F(t,U^sx))=S(F(t+s,x))=G(t+s,x).$$
Now fix $a>0$ and define $F_a(t,x)=F(t,x)$ if $|t|<a$, $F_a(t,x)=0$ otherwise, and let
$$G_a(t,x)=S(F_a(t,x)).$$
Since $S$ is semilocal it has support in $|t|<a+\epsilon$ for a localization constant $\epsilon$ and since $S$ is positive (i.e., its values are non-negative functions) and sublinear, we have
\begin{align*}
G(t,x)=S(F)&=S(F_{a+\epsilon}+(F-F_{a+\epsilon}))\\
& \leq S(F_{a+\epsilon})+S(F-F_{a+\epsilon})
\end{align*}
and since $F-F_{a+\epsilon}$ has support in $|t|>a+\epsilon$, and $S$ is semilocal, the last term on the right vanishes for $|t|\leq a$ for $\epsilon$ sufficiently large, independently of $a$. Thus we have $G\leq G_{a+\epsilon}$ for $|t|\leq a$. Suppose now that there exists a constant $C_1>0$ such that
$$\int_{\mathbb{R}}|Sf(t)|^pw(U^tx)\, dt\leq C_1\int_{\mathbb{R}}|f(t)|^pv(U^tx)\, dt$$
for almost every $x$. Then since $G(0,x)^pw(x)$ and $G(t,x)^pw(U^tx)$ are equimeasurable functions of $x$, we have
\begin{align*}
2\int_XG(0,x)^pw(x)\,d\mu (x)&=\frac{1}{a}\int_{|t|<a}\,dt\int_XG(t,x)^pw(U^tx)\,d\mu (x)\\
&\leq\frac{1}{a}\int_{|t|<a}\,dt\int_XG_{a+\epsilon}(t,x)^pw(U^tx)\,d\mu (x)\\
&=\frac{1}{a}\int_X \int_{|t|<a}G_{a+\epsilon}(t,x)^pw(U^tx)\,dt\, d\mu (x)
\end{align*}
by the Fubini's theorem. Note also that the function $F_{a+\varepsilon}(\cdot,x)$ belongs to $L^1(\mathbb R)$ since $F(t,x)$ is a locally integrable function of $t$ for almost all $x$. Since $SF_{a+\epsilon}=G_{a+\epsilon}$ we have
$$\int_{|t|<a}G_{a+\epsilon}(t,x)^pw(U^tx)\,dt\leq C_1\int_{\mathbb{R}}|F_{a+\epsilon}(t,x)|^pv(U^tx)\,dt.$$
Using the fact that  $F_{a+\epsilon}(t,x)$ is supported in $|t|<a+\epsilon$ and the Fubini's theorem respecively gives
\begin{align*}
2\int_XG(0,x)^pw(x)\,d\mu (x)&\leq\frac{1}{a}C_1\int_X \int_{\mathbb{R}}|F_{a+\epsilon}(t,x) |^pv(U^tx)\,dt\, d\mu (x)\\
&=\frac{1}{a}C_1\int_X \int_{|t|<a+\epsilon}|F_{a+\epsilon}(t,x) |^pv(U^tx)\,dt\, d\mu (x)\\
&=\frac{1}{a}C_1(2a+2\epsilon )\int_X|F_{a+\epsilon}(t,x) |^pv(U^tx)\, d\mu (x)\\
&=\frac{1}{a}C_1(2a+2\epsilon )\int_X|F(t,x) |^pv(U^tx)\, d\mu (x)\\
&=\frac{1}{a}C_1(2a+2\epsilon )\int_X|F(0,x) |^pv(x)\, d\mu (x)
\end{align*}
since $|F(0,x) |^pv(x)$ and $|F(t,x) |^pv(U^tx)$ are equimeasurable. This yields
$$\int_XG(0,x)^pw(x)\,d\mu (x)\leq C_1\int_X|F(0,x) |^pv(x)\, d\mu (x)+\frac{C_1\epsilon}{a}\int_X|F(0,x) |^pv(x)\, d\mu (x).$$
Letting   $a$ tend to infinity gives
$$\int_{X}|S^{\sharp}f(x)|^pw(x)\, d\mu (x)\leq C_1\int_{X}|f(x)|^pv(x)\, d\mu (x),$$
proving the first part of our theorem.\\
\indent
Suppose now that $S$ satisfies a weighted weak type inequality. For any given $\lambda >0$, let $E$ and $\widetilde{E}$ be the sets of points where $ G(0,x)>\lambda$ and $G_{a+\epsilon}(t,x)>\lambda$, respectively, and $\widetilde{E}_y$ the intersection of $\widetilde{E}$ with the set $\{(t,x):x=y\}$. Then we have
$$2aw(E)\leq w(\widetilde{E})=\int_X w(\widetilde{E}_x)\,d\mu (x).$$
On the other hand, $S$ satisfies  weighted weak type inequality
$$w(\widetilde{E}_x)\leq\frac{C_2}{\lambda^p}\int_{\mathbb{R}}|F_{a+\epsilon}(t,x) |^pv(U^tx)\,dt$$
and since $F_{a+\epsilon}(t,x)$ is supported in $|t|<a+\epsilon$, we obatin
\begin{align*}
w(\widetilde{E}_x)&\leq\frac{C_2}{\lambda^p}\int_{\mathbb{R}}|F_{a+\epsilon}(t,x) |^pv(U^tx)\,dt\\
&=\frac{C_2}{\lambda^p}\int_{|t|<a+\epsilon}|F_{a+\epsilon}(t,x) |^pv(U^tx)\,dt\\
\end{align*}
Hence we have
\begin{align*}
2aw(E)&\leq \int_X w(\widetilde{E}_x)\,d\mu (x)\\
&\leq \frac{C_2}{\lambda^p}\int_X\int_{|t|<a+\epsilon}|F_{a+\epsilon}(t,x) |^pv(U^tx)\,dt\, d\mu (x)\\
&= \frac{C_2}{\lambda^p}(2a+2\epsilon )\int_X|F_{a+\epsilon}(t,x) |^pv(U^tx)\,dt\, d\mu (x)\\
&= \frac{C_2}{\lambda^p}(2a+2\epsilon )\int_X|F(t,x) |^pv(U^tx)\,dt\, d\mu (x)\\
&= \frac{C_2}{\lambda^p}(2a+2\epsilon )\int_X|F(0,x) |^pv(x)\,dt\, d\mu (x)\\
\end{align*}
since $|F(0,x) |^pv(x)$ and $|F(t,x) |^pv(U^tx)$ are equimeasurable. Hence we obtain
$$w(E)\leq \frac{C_2}{\lambda^p}\int_X|F(0,x) |^pv(U^tx)\,d\mu (x)+\frac{C_2\epsilon}{a\lambda^p}\int_X|F(0,x) |^pv(x)\,dt\, d\mu (x).$$
When we let $a$ tend to $\infty$ we find the desired result since.
\end{proof}
\indent
The following theorem is our next result:
\begin{thm}\label{wifep2}Let $w$ and $v$ be two non-negative measurable functions and suppose also that $w(U^tx)$ and $v(U^tx)$ are measurable in the product of $X$ with the real line, and let $(S_n)$ and $(P_n)$  be sequences of operators as above. Define
$Sf=\sup_n|S_nf|$, $S^\sharp f=\sup_n|S_n^\sharp f|,$ and $Pf=\sup_n|P_nf|$, $P^\sharp f=\sup_n|P_n^\sharp f|$.  Suppose that there exists a constant $C_1>0$ such that
$$\int_{\mathbb{R}}|Sf(t)|^pw(U^tx)\, dt\leq C_1\int_{\mathbb{R}}|Pf(t)|^pv(U^tx)\, dt$$
and $Pf(U^{(\cdot )}x)\in L^p_{v(U^{(\cdot )}x)}(\mathbb{R})$ for almost every $x$  , for all $f\in (L^1+L^\infty )(\mathbb{R})$,  $1\leq p<\infty$. Then we have
$$\int_{X}|S^{\sharp}f(x)|^pw(x)\, d\mu (x)\leq C_1\int_{X}|P^{\sharp}f(x)|^pv(x)\, d\mu (x) $$
for all  $f\in (L^1+L^\infty )(X)$ satisfying $P^{\sharp}f\in L^p_v(X)$,  $1\leq p<\infty$.\\
\indent
Also, if there exists a constant $C_2>0$ such that for every $\lambda>0$
$$\int_{\{t:|Sf(t)|>\lambda\}}w(U^tx)\, dt\leq \frac{C_2}{\lambda^p}\int_{\mathbb{R}}|Pf(t)|^pv(U^tx)\, dt$$
and $Pf(U^{(\cdot )}x)\in L^p_{v(U^{(\cdot )}x)}(\mathbb{R})$ for almost every $x$  , for all $f\in (L^1+L^\infty )(\mathbb{R})$, then for every $\lambda>0$
$$\int_{\{x:|S^{\sharp}f(x)|>\lambda\}}w(x)\, d\mu (x)\leq \frac{C_2}{\lambda^p}\int_{X}|P^{\sharp}f(x)|^pv(x)\, d\mu (x)$$
holds for all   $f\in (L^1+L^\infty )(X)$ satisfying $P^{\sharp}f\in L^p_v(X)$, $1\leq p<\infty$.
\end{thm}
\begin{proof}We prove this result similar to the proof of Theorem~\ref{wifep} but note that we are transferring two different operators in this case and therefore, we need to build the same type of setup for both operators involved in the statement of our theorem.  More precisely, let
$$G(t,x)=S(F(t,x))$$
and
$$H(t,x)=P(F(t,x)).$$
Then we have $S^{\sharp}f=g(x)$ with $g(x)=G(0,x)$ and $P^{\sharp}f=h(x)$ with $h(x)=H(0,x)$ as in setup given in the introduction.\\
\indent
We have
$$F(t,U^sx)=F(t+s,x),$$
which means that for any two given values $t_1$, $t_2$ of $t$, $F(t_1,x)$ and $F(t_2,x)$ are equimeasurable functions of $x$. Also, because of the  translation invariance of $S$, the function $G(t,x)$ has the same property. Indeed we have
$$G(t,U^sx)=S(F(t,U^sx))=S(F(t+s,x))=G(t+s,x).$$
Now fix $a>0$ and define
$$
F_a(t,x)=
\begin{cases}
F(t,x), & |t|<a,\\
0,& |t|\ge a,
\end{cases}
$$
and let
$$G_a(t,x)=S(F_a(t,x)).$$
Since $S$ and $P$ are semilocal there exist localization constants $\epsilon_s$ and $\epsilon_p$ such that $S$ has support in $|t|<a+\epsilon_s$ and $P$ has support in $|t|<a+\epsilon_p$. Define  $\epsilon =\max (\epsilon_s ,\epsilon_p )$, then both $S$ and $P$ have support in $|t|<a+\epsilon$.\\
\indent
Because $S$ is positive and sublinear, we have
\begin{align*}
G(t,x)=S(F)&=S(F_{a+\epsilon}+(F-F_{a+\epsilon}))\\
& \leq S(F_{a+\epsilon})+S(F-F_{a+\epsilon})
\end{align*}
and since $F-F_{a+\epsilon}$ has support in $|t|>a+\epsilon$, and $S$ is semilocal, the last term on the right vanishes for $|t|\leq a$ for $\epsilon$ sufficiently large, independently of $a$. Thus we have $G\leq G_{a+\epsilon}$ for $|t|\leq a$. \\
\indent
Note that the same argument applies to $H(t,x)$ and $P$, i.e., $G(t,x)$ and $S$ can be replaced with $H(t,x)$ and $P$ respectively in the properties outlined above.\\
\indent
Suppose now that there exists a constant $C_1>0$ such that
$$\int_{\mathbb{R}}|Sf(t)|^pw(U^tx)\, dt\leq C_1\int_{\mathbb{R}}|Pf(t)|^pv(U^tx)\, dt.$$
Then since $G(0,x)^pw(x)$ and $G(t,x)^pw(U^tx)$ are equimeasurable functions of $x$, we have
\begin{align*}
2\int_XG(0,x)^pw(x)\,d\mu (x)&=\frac{1}{a}\int_{|t|<a}\,dt\int_XG(t,x)^pw(U^tx)\,d\mu (x)\\
&\leq\frac{1}{a}\int_{|t|<a}\int_XG_{a+\epsilon}(t,x)^pw(U^tx)\,dt\, d\mu (x)\\
&=\frac{1}{a}\int_X \int_{|t|<a}G_{a+\epsilon}(t,x)^pw(U^tx)\,dt\, d\mu (x).
\end{align*}
For almost all $x$, the function $F_{a+\varepsilon}(\cdot,x)$ belongs to $L^1(\mathbb R)$ since $F(t,x)$ is a locally integrable function of $t$ for almost all $x$. Since $SF_{a+\epsilon}=G_{a+\epsilon}$ and  $PF_{a+\epsilon}=H_{a+\epsilon}$ we have
\begin{align*}
\int_{|t|<a}G_{a+\epsilon}(t,x)^pw(U^tx)\,dt&\leq \int_{\mathbb{R}}G_{a+\epsilon}(t,x)^pw(U^tx)\,dt\\
&\leq C_1\int_{\mathbb{R}}|H_{a+\epsilon}(t,x)|^pv(U^tx)\,dt.
\end{align*}
Since $H_{a+\epsilon}(t,x)$ is supported in $|t|<a+\epsilon$  we obtain
\begin{align*}
2\int_XG(0,x)^pw(x)\,d\mu (x)&\leq\frac{1}{a}C_1\int_X\int_{\mathbb{R}}|H_{a+\epsilon}(t,x)|^pv(U^tx)\,dt\, d\mu (x)\\
&=\frac{1}{a}C_1\int_X \int_{|t|<a+\epsilon}|H_{a+\epsilon}(t,x) |^pv(U^tx)\,dt\, d\mu (x)\\
&=\frac{1}{a}C_1\int_X \int_{|t|<a+\epsilon}|H_{a+\epsilon}(t,x) |^pv(U^tx)\,dt\, d\mu (x)\\
&=\frac{1}{a}C_1\int_{|t|<a+\epsilon}\int_X|H_{a+\epsilon}(t,x) |^pv(U^tx)\, d\mu (x)\,dt\\
&=\frac{1}{a}C_1(2a+2\epsilon )\int_X|H_{a+\epsilon}(t,x) |^pv(U^tx)\, d\mu (x)\\
&=\frac{1}{a}C_1(2a+2\epsilon )\int_X|H(t,x) |^pv(U^tx)\, d\mu (x)\\
&=\frac{1}{a}C_1(2a+2\epsilon )\int_X|H(0,x) |^pv(x)\, d\mu (x)
\end{align*}
since $|H(0,x) |^pv(x)$ and $|H(t,x) |^pv(U^tx)$ are equimeasurable. Thus we obtain
$$\int_XG(0,x)^pw(x)\,d\mu (x)\leq C_1\int_X|H(0,x) |^pv(x)\, d\mu (x)+\frac{C_1\epsilon}{a}\int_X|H(0,x) |^pv(x)\, d\mu (x).$$
Letting $a\to\infty$ gives
$$\int_X |S^\sharp f(x)|^p w(x)\,d\mu(x)\leq C_1 \int_X |P^\sharp f(x)|^p v(x)\,d\mu(x),$$
which proves the first part of our theorem.\\
\indent
Suppose now that $S$ satisfies a weighted weak type inequality. For any given $\lambda >0$, let $E$ and $\widetilde{E}$ be the sets of points where $ G(0,x)>\lambda$ and $G_{a+\epsilon}(t,x)>\lambda$, respectively, and $\widetilde{E}_y$ the intersection of $\widetilde{E}$ with the set $\{(t,x):x=y\}$. Then we have
$$2aw(E)\leq w(\widetilde{E})=\int_X w(\widetilde{E}_x)\,d\mu (x).$$
On the other hand, $S$ satisfies  weighted weak type inequality
$$w(\widetilde{E}_x)\leq\frac{C_2}{\lambda^p}\int_{\mathbb{R}}|H_{a+\epsilon}(t,x) |^pv(U^tx)\,dt.$$
Thus we have
\begin{align*}
2aw(E)&\leq \int_X w(\widetilde{E}_x)\,d\mu (x)\\
&\leq \frac{C_2}{\lambda^p}\int_X \int_{\mathbb{R}}|H_{a+\epsilon}(t,x) |^pv(U^tx)\,dt\, d\mu (x)\\
&= \frac{C_2}{\lambda^p}\int_X \int_{|t|<a+\epsilon}|H_{a+\epsilon}(t,x) |^pv(U^tx)\,dt\, d\mu (x)\\
&= \frac{C_2}{\lambda^p}\int_{|t|<a+\epsilon}\int_X|H_{a+\epsilon}(t,x) |^pv(U^tx)\, d\mu (x)\,dt\\
&= \frac{C_2}{\lambda^p}(2a+2\epsilon )\int_X|H_{a+\epsilon}(t,x) |^pv(U^tx)\, d\mu (x)\\
&= \frac{C_2}{\lambda^p}(2a+2\epsilon )\int_X|H(t,x) |^pv(U^tx)\, d\mu (x)\\
&= \frac{C_2}{\lambda^p}(2a+2\epsilon )\int_X|H(0,x) |^pv(x)\, d\mu (x)
\end{align*}
since $|H(0,x) |^pv(x)$ and $|H(t,x) |^pv(U^tx)$ are equimeasurable.  This gives
\begin{align*}
w(E)\leq \frac{C_2}{\lambda^p}\int_X|H(0,x) |^pv(U^tx)\,d\mu (x)+\frac{C_2\epsilon}{a\lambda^p}\int_X|H(0,x) |^pv(U^tx)\,d\mu (x)
\end{align*}
yielding
$$\int_{\{x:|S^{\sharp}f(x)|>\lambda\}}w(x)\, d\mu\leq \frac{C_2}{\lambda^p}\int_{X}|P^{\sharp}f(x)|^pv(x)\, d\mu$$
when we let $a$ tend to $\infty$ which completes the proof of our theorem.
\end{proof}
\indent
Our next result is the following vector-valued extension of Theorem~\ref{wifep}:
 \begin{thm}\label{vvctpws}Let $w$ and $v$ be two non-negative measurable functions and suppose also that $w(U^tx)$ and $v(U^tx)$ are measurable in the product of $X$ with the real line, and let $1\leq r\leq\infty$ and $1\leq p<\infty$. Also, let  $T_n$ be a sequence of operators as described in the introduction and $Sf=\sup|T_nf|$.  Suppose that there exists a constant $C_1>0$ such that
 $$\left(\int_{\mathbb{R}}\|Sf(t)\|_{\ell^r}^pw(U^tx)\,dt\right)^{1/p}\leq C_1\left(\int_{\mathbb{R}}\|f(t)\|_{\ell^r}^pv(U^tx)\,dt\right)^{1/p}$$
 for almost every $x$, where $f=(f_1,f_2,f_3,\dots )$ is a sequence of functions in $(L^1+L^\infty )(\mathbb{R})$ such that $\|f\|_{\ell^r}\in (L^1+L^\infty )(\mathbb{R})$.
 \indent
Then we have
  $$\left(\int_X\|S^\sharp f(x)\|_{\ell^r}^pw(x)\,d\mu (x)\right)^{1/p}\leq C_1\left(\int_X\|f(x)\|_{\ell^r}^pv(x)\,d\mu (x)\right)^{1/p}$$
  where $f=(f_1,f_2,f_3,\dots )$ is a sequence of functions in $(L^1+L^\infty )(X)$ such that $f\in L_{\ell^r}^p(v)$.
\indent
 Also, if there exists a constant $C_2>0$ such that for every $\lambda>0$
$$\int_{\{t:\|Sf(t)\|_{\ell^r}>\lambda\}}w(U^tx)\, dt\leq \frac{C_2}{\lambda^p}\int_{\mathbb{R}}\|f(t)\|_{\ell^r}^pv(U^tx)\, dt$$
holds for almost every $x$, $1\leq p<\infty$, and for all sequences of functions $f=(f_1,f_2,f_3,\dots )$ in $(L^1+L^\infty )$ such that $\|f\|_{\ell^r}\in (L^1+L^\infty )(\mathbb{R})$, then for every $\lambda>0$
$$\int_{\{x:\|S^\sharp f(x)\|_{\ell^r}>\lambda\}}w(x)\, d\mu (x)\leq \frac{C_2}{\lambda^p}\int_{X}\|f(x)\|_{\ell^r}^pv(x)\, d\mu (x)$$
holds for all sequences of functions $f=(f_1,f_2,f_3,\dots )$ in $ (L^1+L^\infty )(X)$ such that $f\in L_{\ell^r}^p(v)$, $1\leq p<\infty$.
\end{thm}
\begin{proof} Under this assumption the operator $S$ has the same properties as the operator $T$ in the setup and as in the proof of Theorem~\ref{wifep} for each term $f_i$ of the sequence $f=(f_1,f_2,f_3,\dots )$ we have
$$F^i(t,U^sx)=F^i(t+s,x),$$
which means that for any two given values $t_1$, $t_2$ of $t$, $F^i(t_1,x)$ and $F^i(t_2,x)$ are equimeasurable functions of $x$. On the other hand, due to translation invariance of $S$, the function $G^i(t,x)$ has the same property. In fact we have
$$G^i(t,U^sx)=S(F^i(t,U^sx))=S(F^i(t+s,x))=G^i(t+s,x).$$
Now fix $a>0$ and define $F^i_a(t,x)=F^i(t,x)$ if $|t|<a$, $F^i_a(t,x)=0$ otherwise, and let
$$G^i_a(t,x)=S(F^i_a(t,x)).$$
\indent
Since $S$ is semilocal the support of $S$ is always contained in an $\epsilon$-neighborhood of the support of the function $\|f\|_{\ell^r}$ by the definition. Let this neighborhood be $|t|>a+\epsilon$ for a fixed $a$ and $\epsilon$. Because of the definition of $\ell^r$ norm we see that each term $f_i$ of the sequence also has support in this neighborhood otherwise the summation of the absolute values does not vanish. To see this assume that for $t_0\in\mathbb{R}$ we have $\|f(t_0)\|_{\ell^r}=0$ and $|f_j(t_0)|>0$ for some $j$. Then we have
$$0<|f_j(t_0)|^r\leq \sum_{k=1}^\infty |f_k(t_0)|^r=0$$
which is impossible.\\
\indent
Now define 
$$
F^i_a(t,x)=
\begin{cases}
F^i(t,x), & |t|<a,\\
0,& |t|\ge a,
\end{cases}
$$
for each $i=1,2,3,\dots$ and let $G^i_a(t,x)=S(F^i_a(t,x))$. Since $S$ is positive (i.e., its values are non-negative functions) and sublinear, we have
\begin{align*}
G^i(t,x)=S(F^i)&=S(F^i_{a+\epsilon}+(F^i-F^i_{a+\epsilon}))\\
& \leq S(F^i_{a+\epsilon})+S(F^i-F^i_{a+\epsilon})
\end{align*}
and since $F^i-F^i_{a+\epsilon}$ has support in $|t|>a+\epsilon$, and $S$ is semilocal, the last term on the right vanishes for $|t|\leq a$ for $\epsilon$ sufficiently large, independently of $a$. Thus we have $G^i\leq G^i_{a+\epsilon}$ for $|t|\leq a$. Suppose now that there exists a constant $C_1>0$ such that
$$\int_{\mathbb{R}}\|Sf(t)\|_{\ell^r}^pw(U^tx)\, dt\leq C_1\int_{\mathbb{R}}\|f(t)\|_{\ell^r}^pv(U^tx)\, dt$$
for almost every $x$. \\
\indent
Since $G^i(0,x)$ and $G^i(t,x)$ are equimeasurable functions of $x$, it is not hard to see that $\|G(0,x)\|^p_{\ell^r}w(x)$ and $\|G(t,x)\|^p_{\ell^r}w(U^tx)$ are also equimeasurable functions of $x$, where 
$$G(t,x)=(G^1(t,x),G^2(t,x),G^3(t,x),\dots ).$$
Thus we obtain
\begin{align*}
2\int_X\|G(0,x)\|_{\ell^r}^pw(x)\,d\mu (x)&=\frac{1}{a}\int_{|t|<a}\,dt\int_X\|G(t,x)\|_{\ell^r}^pw(U^tx)\,d\mu (x)\\
&\leq\frac{1}{a}\int_{|t|<a}\,dt\int_X\|G_{a+\epsilon}(t,x)\|_{\ell^r}^pw(U^tx)\,d\mu (x)\\
&=\frac{1}{a}\int_X \int_{|t|<a}\|G_{a+\epsilon}(t,x)\|_{\ell^r}^pw(U^tx)\,dt\, d\mu (x).
\end{align*}
\indent
Also note that for each $i=1,2,3,\dots $ the function $F^i_{a+\varepsilon}(\cdot,x)$ belongs to $L^1(\mathbb R)$ since $F^i(t,x)$ is a locally integrable function of $t$ for almost all $x$ and $\|F_{a+\varepsilon}(\cdot,x)\|_{\ell^r}$ belongs to $L^1(\mathbb R)$ since  $\|F(t,x)\|_{\ell^r}$   is a locally integrable function of $t$ for almost all $x$, where  
$$F(t,x)=(F^1(t,x), F^2(t,x),F^3(t,x),\dots ).$$
\indent
Since $SF^i_{a+\epsilon}=G^i_{a+\epsilon}$ we have
$$\int_{|t|<a}\|G_{a+\epsilon}(t,x)\|_{\ell^r}^pw(U^tx)\,dt\leq C_1\int_{\mathbb{R}}\|F_{a+\epsilon}(t,x)\|_{\ell^r}^pv(U^tx)\,dt,$$
and since $\|F_{a+\epsilon}(t,x)\|_{\ell^r}$ is supported in $|t|<a+\epsilon$, we obtain
\begin{align*}
2\int_X\|G(0,x)\|_{\ell^r}^pw(x)\,d\mu (x)&\leq\frac{1}{a}C_1\int_X \int_{\mathbb{R}}\|F_{a+\epsilon}(t,x)\|_{\ell^r}^pv(U^tx)\,dt\, d\mu (x)\\
&=\frac{1}{a}C_1\int_X \int_{|t|<a+\epsilon}\|F_{a+\epsilon}(t,x)\|_{\ell^r}^pv(U^tx)\,dt\, d\mu (x)\\
&=\frac{1}{a}C_1\int_{|t|<a+\epsilon}\int_X \|F_{a+\epsilon}(t,x)\|_{\ell^r}^pv(U^tx)\,d\mu (x)\, dt\\
&=\frac{1}{a}C_1(2a+2\epsilon )\int_X \|F_{a+\epsilon}(t,x)\|_{\ell^r}^pv(U^tx)\,d\mu (x)\\
&=\frac{1}{a}C_1(2a+2\epsilon )\int_X \|F(t,x)\|_{\ell^r}^pv(U^tx)\,d\mu (x)\\
&=\frac{1}{a}C_1(2a+2\epsilon )\int_X \|F(0,x)\|_{\ell^r}^pv(x)\,d\mu (x)
\end{align*}
since $\|F(0,x)\|_{\ell^r}^pv(x)$ and $\|F(t,x)\|_{\ell^r}^pv(U^tx)$ are equimeasurable. This yields
$$\int_X\|G(0,x)\|_{\ell^r}^pw(x)\,d\mu (x)\leq C_1\int_X \|F(0,x)\|_{\ell^r}^pv(x)\,d\mu (x)+\frac{C_1\epsilon}{a}\int_X \|F(0,x)\|_{\ell^r}^pv(x)\,d\mu (x).$$
Letting  $a$ tend to infinity yields
$$\int_{X}\|S^{\sharp}f(x)\|_{\ell^r}^pw(x)\, d\mu (x)\leq C_1\int_{X}\|f(x)\|_{\ell^r}^pv(x)\, d\mu (x),$$
proving the first part of our theorem for $1\leq r <\infty$.\\
\indent
For $r=\infty$, we have  $\|f\|_{\ell^\infty}=\sup_i|f_i|$. Again using the semilocallity of $S$ we see that the support of $S$ is contained in a neighborhood $|t|>a+\epsilon$ for a fixed $a$ and an $\epsilon$ of the support of the function $\|f\|_{\ell^\infty}$. This clearly means that for each $i=1,2,3,\dots$ the function $f_i$ is supported in this neighborhood as well,  since we take the supremum of absolute values, i.e., when the supremum  vanishes at any point all terms have to vanish at that point. The rest of the proof for the case $r=\infty$ can be completed as way we have proved the theorem for the case $1\leq r<\infty$.\\
\indent
Suppose now that $S$ satisfies a weighted weak type inequality. For any given $\lambda >0$, let $E$ and $\widetilde{E}$ be the sets of points where $\|G(0,x)\|_{\ell^r}>\lambda$ and $\|G_{a+\epsilon}(t,x)\|_{\ell^r}>\lambda$, respectively, and $\widetilde{E}_y$ the intersection of $\widetilde{E}$ with the set $\{(t,x):x=y\}$. Then we have
$$2aw(E)\leq w(\widetilde{E})=\int_X w(\widetilde{E}_x)\,d\mu (x).$$
On the other hand, $S$ satisfies  weighted weak type inequality
$$w(\widetilde{E}_x)\leq\frac{C_2}{\lambda^p}\int_{\mathbb{R}}\|F_{a+\epsilon}(t,x)\|_{\ell^r}^pv(U^tx)\,dt$$
and since $F_{a+\epsilon}(t,x)$ is supported in $|t|<a+\epsilon$, we obtain
\begin{align*}
w(\widetilde{E}_x)&\leq\frac{C_2}{\lambda^p}\int_{\mathbb{R}}\|F_{a+\epsilon}(t,x)\|_{\ell^r}^pv(U^tx)\,dt\\
&=\frac{C_2}{\lambda^p}\int_{|t|<a+\epsilon}\|F_{a+\epsilon}(t,x)\|_{\ell^r}^pv(U^tx)\,dt.
\end{align*}
Hence we have
\begin{align*}
2aw(E)&\leq \int_X w(\widetilde{E}_x)\,d\mu (x)\\
&\leq \frac{C_2}{\lambda^p}\int_X\int_{|t|<a+\epsilon}\|F_{a+\epsilon}(t,x)\|_{\ell^r}^pv(U^tx)\,dt\, d\mu (x)\\
&= \frac{C_2}{\lambda^p}\int_{|t|<a+\epsilon}\int_X\|F_{a+\epsilon}(t,x)\|_{\ell^r}^pv(U^tx)\,d\mu (x)\,dt\\
&= \frac{C_2}{\lambda^p}(2a+2\epsilon )\int_X\|F_{a+\epsilon}(t,x)\|_{\ell^r}^pv(U^tx)\,d\mu (x)\,dt\\
&= \frac{C_2}{\lambda^p}(2a+2\epsilon )\int_X\|F(t,x)\|_{\ell^r}^pv(U^tx)\,d\mu (x)\,dt\\
&= \frac{C_22a}{\lambda^p}\int_X\|F(0,x)\|_{\ell^r}^pv(x)\,d\mu (x)+\frac{C_22\epsilon}{\lambda^p}\int_X\|F(0,x)\|_{\ell^r}^pv(x)\,d\mu (x)
\end{align*}
since  $\|F(0,x)\|_{\ell^r}^pv(x)$ and $\|F(t,x)\|_{\ell^r}^pv(U^tx)$ are equimeasurable. Hence we have
$$w(E)\leq \frac{C_2}{\lambda^p}\int_X\|F(0,x)\|_{\ell^r}^pv(x)\,d\mu (x)+\frac{C_2\epsilon}{a\lambda^p}\int_X\|F(0,x)\|_{\ell^r}^pv(x)\,d\mu (x).$$
When we let $a$ tend to $\infty$ we find the desired result for the case $1\leq r <\infty$.\\
\indent
The case $r=\infty$ is being proved as in the first part of our proof.
\end{proof}
\indent
We also have the following extension of Theorem~\ref{wifep2} to the vector-valued setting:
\begin{thm}\label{vvwifep6}Let $w$ and $v$ be two non-negative measurable functions and suppose also that $w(U^tx)$ and $v(U^tx)$ are measurable in the product of $X$ with the real line, and let $(S_n)$ and $(P_n)$  be sequences of operators as above. Define
$Sf=\sup_n|S_nf|$, $S^\sharp f=\sup_n|S_n^\sharp f|,$ and $Pf=\sup_n|P_nf|$, $P^\sharp f=\sup_n|P_n^\sharp f|$.  Suppose that for all $f\in (L^1+L^\infty )(\mathbb{R})$,  $\|Pf(U^tx)\|_{\ell^r}$ is a locally integrable function of $t$ for almost all $x$ and that there exists a constant $C_1>0$ such that
$$\int_{\mathbb{R}}\|Sf(t)\|_{\ell^r}^pw(U^tx)\, dt\leq C_1\int_{\mathbb{R}}\|Pf(t)\|_{\ell^r}^pv(U^tx)\, dt$$
for almost every $x$, for all sequences  $f=(f_1,f_2,f_3,\dots )$ in  $(L^1+L^\infty )(\mathbb{R})$ such that $\|f\|_{\ell^r}\in (L^1+L^\infty )(\mathbb{R})$, and $1\leq p<\infty$. Then we have
$$\int_{X}\|S^{\sharp}f(x)\|_{\ell^r}^pw(x)\, d\mu (x)\leq C_1\int_{X}\|P^{\sharp}f(x)\|_{\ell^r}^pv(x)\, d\mu (x) $$
for all sequences  $f=(f_1,f_2,f_3,\dots )$ in  $(L^1+L^\infty )(X)$ such that $P^\sharp f\in L_{\ell^r}^p(v)$.\\
\indent
Also, if there exists a constant $C_2>0$ such that for every $\lambda>0$
$$\int_{\{t:\|Sf(t)\|_{\ell^r}>\lambda\}}w(U^tx)\, dt\leq \frac{C_2}{\lambda^p}\int_{\mathbb{R}}\|Pf(t)\|_{\ell^r}^pv(U^tx)\, dt$$
holds for almost every $x$, $1\leq p<\infty$, and for all sequences of functions $f=(f_1,f_2,f_3,\dots )$ in $(L^1+L^\infty )(\mathbb{R})$ such that $\|f\|_{\ell^r}\in (L^1+L^\infty )(\mathbb{R})$, then for every $1\leq p<\infty$ and  $\lambda>0$
$$\int_{\{x:\|S^{\sharp}f(x)\|_{\ell^r}>\lambda\}}w(x)\, d\mu (x)\leq \frac{C_2}{\lambda^p}\int_{X}\|P^{\sharp}f(x)\|_{\ell^r}^pv(x)\, d\mu (x)$$
holds for all   $f\in (L^1+L^\infty )(X)$ such that $P^\sharp f\in L_{\ell^r}^p(v)$.
\end{thm}
\begin{proof}This result can be proved by a close inspection of the arguments we used in the proofs Theorem~\ref{wifep2} and Theorem~\ref{vvctpws}.
\end{proof}
\section{Applications to Ergodic Theory}
\indent

In this section we apply the results of the previous section to some weighted inequalities on the real line also proved by the author to prove some new inequalities in ergodic theory.\\
\indent
We say that a non-negative measurable function $w\in L_{{\rm{loc}}}^1(\mathbb{R})$ is an $A_p$ weight for some $1<p<\infty$ if
\begin{align*}
\sup_I\left[\frac{1}{|I|}\int_Iw(x)\,dx\right]\left[\frac{1}{|I|}\int_I[w(x)]^{-1/(p-1)}\,dx\right]^{p-1}\leq C<\infty .
\end{align*}
The supremum is taken over all intervals $I\subset \mathbb{R}$; $|I|$ denotes the measure of $I$.\\
\indent
We say that $w\in A_1$ if
$$\frac{1}{|I|}\int_I w(x)\, dx\leq C\,{\rm{ess\, inf}}_I\,w$$
for all intervals $I$.\\
\indent
Note that the class $A_\infty$ is defined as
 $$A_{\infty}=\cup_{p\geq 1}A_p.$$
Thus, $w\in A_\infty$ implies $w\in A_p$  for some  $1\leq p<\infty$, and if $w\in A_p$ for any $1\leq p<\infty$, then $w\in A_{\infty}$. If $w$ is an $A_\infty$ weight then there exist $\delta >0$ and $\epsilon >0$ such that, given an interval $I\subset \mathbb{R}$ and 
 for any measurable set $E\subset I$, we have
 $$\frac{w(E)}{w(I)}\leq \delta \left(\frac{|E|}{|I|}\right)^{\epsilon},$$
 where
 $$w(A)=\int_Aw(x)\,dx$$
 for any subset $A\subset \mathbb{R}^n$.\\
\indent
Also, if $w\in A_p$ then there exist $u, v\in A_1$ such that $w=uv^{1-p}$.\\
Before presenting our main result  let us first mention some of the well-known properties of the functions in $A_p^{\lambda}$.\\
\indent
Let $\lambda$ be a real number with $\lambda\geq 1$,\\
\indent
(i)  Given $p$, $1<p<\lambda^\prime$, we  say that a weight $w$ belongs to $A_p^{\lambda}$ if there exists a constant $C$ such that
$$\sup_I\left(\frac{1}{|I|}\int_Iw(x)^{\lambda^\prime/(\lambda^\prime-p)}\, dx\right)^{(\lambda^\prime-p)/\lambda^\prime}\left(\frac{1}{|I|}\int_Iw(x)^{-1/(p-1)}\, dx\right)^{p-1}\leq C$$
where $I$ stands for an arbitrary interval in $\mathbb{R}$.\\
\indent
(ii) We say that $w\in A_1^{\lambda}$ if there is a constant $C$ such that
$$M(w^\lambda )(x)\leq Cw^\lambda (x),\;\;\; {\textrm{a.e.}}\; x,$$
where $M$ is the Hardy-Littlewood maximal function defined as
$$Mf(x)=\sup_{x\in I}\frac{1}{|I|}\int_I|f(y)|\, dy$$
where $I$ stands for an arbitrary interval in $\mathbb{R}$.\\
\indent
Note also that the class $A_p^1$ coincides with Muckenhoupt's class $A_p$, and the following duality property holds:
$$w\in A_p^\lambda \Leftrightarrow w^{-p^\prime /p}\in A_{p^\prime /\lambda}, \;\;\;\; \lambda \geq1,\;\; 1<p<\lambda^\prime .$$
\indent
The readers are referred to J.~Garcia-Cuerva and J. L.~Rubio de Francia~\cite{jgcjrdf} for more information on the classes $A_p$,
 and to D. S.~Kurtz and R. L.~Wheeden~\cite{dskrlw} on the classes $A_p^\lambda$.\\
\indent
Let now $(X,\mathscr{B} ,\mu )$ be a totally  $\sigma$-finite measure space, and $\{U_t:-\infty <t<\infty\}$ be a one-parameter group of measure-preserving transformations on $X$.\\
\indent
We say that $w$ satisfies condition $A_p^\prime$ if there exists a constant $M>0$ such that for a.e. $x\in X$ 
\begin{align*}
\sup_I\left[\frac{1}{|I|}\int_Iw(U_tx)\,dt\right]\left[\frac{1}{|I|}\int_I[w(U_tx)]^{-1/(p-1)}\,dt\right]^{p-1}\leq M.
\end{align*}
The supremum is taken over all intervals $I\subset \mathbb{R}$; $|I|$ denotes the measure of $I$.\\
Here 
$w\in A_1^\prime$ if
$$\frac{1}{|I|}\int_I w(U_tx)\, dt\leq C\,{\rm{ess\, inf}}_I\,w(U_tx),\;\;\; {\textrm{a.e.}}\; x,$$
for all intervals $I$.\\
\indent
Similar to the classical case we can also define the class $(A_p^{\lambda})^\prime$ in ergodic theory as follows:\\
\indent
Let $\lambda$ be a real number with $\lambda\geq 1$,\\
\indent
(i)  Given $p$, $1<p<\lambda^\prime$, we  say that a weight $w$ belongs to $(A_p^{\lambda})^\prime$ if there exists a constant $C$ such that for a.e. $x$
$$\sup_I\left(\frac{1}{|I|}\int_Iw(U_tx)^{\lambda^\prime/(\lambda^\prime-p)}\, dt\right)^{(\lambda^\prime-p)/\lambda^\prime}\left(\frac{1}{|I|}\int_Iw(U_tx)^{-1/(p-1)}\, dt\right)^{p-1}\leq C$$
where $I$ stands for an arbitrary interval in $\mathbb{R}$.\\
\indent
(ii) We say that $w\in (A_1^{\lambda})^\prime$ if there is a constant $C$ such that
$$M^\ast (w^\lambda )(x)\leq Cw^\lambda (x),\;\;\; {\textrm{a.e.}}\; x,$$
where $M^\ast$ is the ergodic maximal function defined as
$$M^\ast f(x)=\sup_{x\in I}\frac{1}{|I|}\int_I|f(U_tx)|\, dt$$
where $I$ stands for an arbitrary interval in $\mathbb{R}$.\\
\indent
It is known (see \cite{kp1}) that if $f\in L^1$ and $\lambda >0$, then
$$\mu\{M^\ast f>\lambda\}\leq \frac{1}{\lambda}\int_{\{M^\ast f>\lambda\}}|f|\, d\mu .$$
\indent
We say that $w$ satisfies condition $A_p^\prime$ if there exists a constant $C$ such that for a.e. $x\in X$ 
\indent
Also it is not hard to see as in the classical case that the class $(A_p^1)^\prime$ coincides with the class $A^{\prime}_p$, and the following duality property holds:
$$w\in (A_p^\lambda)^\prime \Leftrightarrow w^{-p^\prime /p}\in A^{\prime}_{p^\prime /\lambda}, \;\;\;\; \lambda \geq1,\;\; 1<p<\lambda^\prime .$$
\indent
Let $f$ be a measurable function defined on $\mathbb{R}$. For each $n\in\mathbb{Z}$ define the operator $A_n$ by
$$A_nf(x)=\frac{1}{2^n}\int_x^{x+2^n}f(y)\, dy.$$
It is a well-known problem to study the different kinds of convergence of the sequence $\{A_nf\}_n$ when $f\in L^p(\mathbb{R})$ for some $1\leq p<\infty$. Consider the variation operator
$$\mathcal{V}f(x)=\left(\sum_{n=-\infty}^\infty|A_nf(x)-A_{n-1}f(x)|^s\right)^{1/s}$$
for $2\leq s<\infty$.\\
\indent
The following result has been obtained by the author in S. Demir~\cite{sdemir}:
\begin{thm}\label{vecapfvar} Let $1<r<\infty$ and $r^\prime =\frac{r}{r-1}$, and suppose that  $2\leq s<\infty$. Then for all $1<\rho <\infty$, and for all $(f_j)\in L^p(w)\cap L^p(\mathbb{R})$, the weighted inequalities
$$\left\|\left(\sum_j(\mathcal{V}f_j)^{\rho}\right)^{1/\rho}\right\|_{L^p(w)}\leq C_{p,\rho}(w)\left\|\left(\sum_j|f_j|^{\rho}\right)^{1/{\rho}}\right\|_{L^p(w)}$$
hold if $w\in A_{p/r^\prime}$ and $r^\prime\leq p<\infty$, or if $w\in A_p^{r^\prime}$ and $1<p< r$. Likewise, if $w(x)^{r^\prime}\in A_1$, then the weak type inequality
\begin{align*}
w\left(\left\{x:\left(\sum_j(\mathcal{V}f_j(x))^{\rho}\right)^{1/\rho}>\lambda\right\}\right)\leq C_{\rho}(w)\frac{1}{\lambda}\int\left(\sum_j|f_j(x)|^{\rho}\right)^{1/\rho}w(x)\, dx
\end{align*}
holds for all $(f_j)\in L^1(w)\cap L^1(\mathbb{R})$.
\end{thm}
\indent
Let now $(X,\mathscr{B} ,\mu )$ be a totally  $\sigma$-finite measure space, and $\{U_t:-\infty <t<\infty\}$ be a one-parameter group of measure-preserving transformations on $X$ as above, and define the ergodic averages
$$B_nf(x)=\frac{1}{n}\int_0^nf(U_{-t}x)\, dt.$$
Let $2\leq s<\infty$,  and define the two-sided variation operator 
$$\mathcal{V}^\prime f(x)=\left(\sum_{n=-\infty}^\infty|B_{2^n}f(x)-B_{2^{n-1}}f(x)|^s\right)^{1/s}$$
for a locally integrable function $f$.\\
\indent
Let $w\in A_p^\prime$ for some $p$, define $U(t,x)=U_t(x)$. It is then not hard to verify that $w\circ U(t,x)\in A_p$ for a.e. $x$. Similarly, if $w\in (A_p^{\lambda})^{\prime}$, then $w\circ U(t,x)\in A_p^{\lambda}$ for a.e. $x$. \\
\indent
With this weight function we can reformulate Theorem~\ref{vecapfvar} as follows:
\begin{thm}\label{vecapfvar2} Let $1<r<\infty$ and $r^\prime =\frac{r}{r-1}$, and suppose that  $2\leq s<\infty$. Then for all $1<\rho <\infty$, and for all $(f_j)\in L^p(w\circ U)\cap L^p(\mathbb{R})$, for a.e. $x$ the weighted inequalities
$$\left(\int \left(\sum_j(\mathcal{V}f_j(t))^{\rho}\right)^{p/\rho}w(U_tx)\, dt\right)^{1/p}\leq C_{p,\rho}(w\circ U)\left(\int\left(\sum_j|f_j(t)|^{\rho}\right)^{p/{\rho}}w(U_tx)\, dt\right)^{1/p}$$
hold if for a.e. x, $w\circ U(t,x) \in A_{p/r^\prime}$ and $r^\prime\leq p<\infty$, or if for a.e. x,  $w\circ U(t,x)\in A_p^{r^\prime}$ and $1<p< r$. Likewise, if for a.e. x,  $w\circ U (t, x)^{r^\prime}\in A_1$, then for a.e. $x$, the weak type inequality
\begin{align*}
\int_{\left\{t:\left(\sum_j(\mathcal{V}f_j(t))^{\rho}\right)^{1/\rho}>\lambda\right\}}w(U_tx)\, dt  \leq C_{\rho}(w\circ U)\frac{1}{\lambda}\int\left(\sum_j|f_j(t)|^{\rho}\right)^{1/\rho}w(U_tx)\, dt
\end{align*}
holds for all $(f_j)\in L^1(w\circ U)\cap L^1(\mathbb{R})$.
\end{thm}
\indent
When we apply Theorem~\ref{vvctpws} to Theorem~\ref{vecapfvar2} we obtain the following result:
\begin{thm}\label{dyvarergo} Let $1<r<\infty$ and $r^\prime =\frac{r}{r-1}$, and suppose that  $2\leq s<\infty$. Then for all $1<\rho <\infty$, and for all $(f_j)\in (L^1+L^\infty )\cap L^p_{\ell^{\rho}}(w)$ , the weighted inequalities
$$\left(\int \left(\sum_j(\mathcal{V}^\prime f_j(x))^{\rho}\right)^{p/\rho}w(x)\, dx\right)^{1/p}\leq C_{p,\rho}(w)\left(\int\left(\sum_j|f_j(x)|^{\rho}\right)^{p/{\rho}}w(x)\, dx\right)^{1/p}$$
hold if $w\in A^{\prime}_{p/r^\prime}$ and $r^\prime\leq p<\infty$, or if $w\in (A_p^{r^\prime})^\prime$ and $1<p< r$. Likewise, if  $w(x)^{r^\prime}\in A^{\prime}_1$, then the weak type inequality
\begin{align*}
\int_{\left\{x:\left(\sum_j(\mathcal{V}^\prime f_j(x))^{\rho}\right)^{1/\rho}>\lambda\right\}}w(x)\, dx  \leq C_{\rho}(w)\frac{1}{\lambda}\int\left(\sum_j|f_j(x)|^{\rho}\right)^{1/\rho}w(x)\, dx
\end{align*}
holds for all $(f_j)\in (L^1+L^\infty )\cap L^1_{\ell^{\rho}}(w)$.
\end{thm}
\indent
In particular, we obtain the following corollary:
\begin{cor}\label{dycorvarergo} Let $1<r<\infty$ and $r^\prime =\frac{r}{r-1}$, and suppose that  $2\leq s<\infty$. Then, for all $f \in (L^1+L^\infty )\cap L_w^p(X)$ the weighted inequalities
$$\left\|\mathcal{V}^\prime f\right\|_{L^p(w)}\leq C_{p}(w)\left\|f\right\|_{L^p(w)}$$
hold if $w \in A^{\prime}_{p/r^\prime}$ and $r^\prime\leq p<\infty$, or if $w\in (A_p^{r^\prime})^\prime$ and $1<p<r$. Likewise, if  $w(x)^{r^\prime}\in A^{\prime}_1$, then the weak type inequality
\begin{align*}
\int_{\left\{x:\mathcal{V}^\prime f(x)>\lambda\right\}}w(x)\, dx  \leq C(w)\frac{1}{\lambda}\left\|f\right\|_{L^1(w)}
\end{align*}
holds for all $f \in (L^1+L^\infty )\cap L_w^1(X)$
\end{cor}
\indent
Recall that a sequence $(n_k)$ of positive integers is called lacunary if there exists a real number $\beta$ such that
$$\frac{n_{k+1}}{n_k}\geq \beta >1$$
for all $k=0, 1,2,3,\dots$\\
\indent
Let now $f$ be a locally integrable function defined on $\mathbb{R}$, and let $(n_k)$ be a lacunary sequence. Define the operator $A_{n_k}$ by
$$A_{n_k}f(x)=\frac{1}{n_k}\int_0^{n_k}f(x-t)\, dt.$$
Define the variation operator
$$\mathcal{V}_sf(x)=\left(\sum_{k=1}^\infty|A_{n_k}f(x)-A_{n_{k-1}}f(x)|^s\right)^{1/s}$$
for $2\leq s<\infty$.\\
\indent
The following result has been obtained by the author in S. Demir~\cite{sdemirvar}:
\begin{thm}\label{vecvvarlrlac}Let $(n_k)$ be a lacunary sequence, $1<r<\infty$, and $r^\prime =\frac{r}{r-1}$. Suppose that  $2\leq s<\infty$. Then, for all $1<\rho <\infty$, and for all $(f_j)\in L^p(w)\cap L^p(\mathbb{R})$, the weighted inequalities
$$\left\|\left(\sum_j(\mathcal{V}_sf_j)^{\rho}\right)^{1/\rho}\right\|_{L^p(w)}\leq C_{p,\rho}(w)\left\|\left(\sum_j|f_j|^{\rho}\right)^{1/{\rho}}\right\|_{L^p(w)}$$
hold if $w\in A_{p/r^\prime}$ and $r^\prime\leq p<\infty$, or if $w\in A_p^{r^\prime}$ and $1<p< r$. Likewise, if $w(x)^{r^\prime}\in A_1$, then the weak type inequality
\begin{align*}
w\left(\left\{x:\left(\sum_j(\mathcal{V}_sf_j(x))^{\rho}\right)^{1/\rho}>\lambda\right\}\right)\leq C_{\rho}(w)\frac{1}{\lambda}\int\left(\sum_j|f_j(x)|^{\rho}\right)^{1/\rho}w(x)\, dx
\end{align*}
holds for all $(f_j)\in L^1(w)\cap L^1(\mathbb{R})$.
\end{thm}
Consider the ergodic averages
$$B_nf(x)=\frac{1}{n}\int_0^nf(U_{-t}x)\, dt$$
as above, and define the variation operator 
$$\mathcal{V}_s^{\prime}f(x)=\left(\sum_{k=1}^\infty|B_{n_k}f(x)-B_{n_{k-1}}f(x)|^s\right)^{1/s}$$
for a locally integrable function $f$.\\
\indent
When we apply Theorem~\ref{vvctpws} to Theorem~\ref{vecvvarlrlac} we obtain the following result:
\begin{thm}\label{vecvvxlac}Let $(n_k)$ be a lacunary sequence,  $1<r<\infty$, and $r^\prime =\frac{r}{r-1}$. Suppose that  $2\leq s<\infty$. Then, for all $1<\rho <\infty$, and for all $(f_j)\in (L^1+L^\infty )\cap L^p_{\ell^{\rho}}(w)$ the weighted inequalities
$$\left(\int \left(\sum_j(\mathcal{V}_s^{\prime} f_j)^{\rho}\right)^{p/\rho}w(x)\, dx\right)^{1/p}\leq C_{p,\rho}(w)\left(\int\left(\sum_j|f_j|^{\rho}\right)^{p/{\rho}}w(x)\, dx\right)^{1/p}$$
hold if $w\in A^{\prime}_{p/r^\prime}$ and $r^\prime\leq p<\infty$, or if $w\in (A_p^{r^\prime})^\prime$ and $1<p<r$. Likewise, if  $w(x)^{r^\prime}\in A^{\prime}_1$, then the weak type inequality
\begin{align*}
\int_{\left\{x:\left(\sum_j(\mathcal{V}_s^{\prime} f_j(x))^{\rho}\right)^{1/\rho}>\lambda\right\}}w(x)\, dx  \leq C_{\rho}(w)\frac{1}{\lambda}\int\left(\sum_j|f_j(x)|^{\rho}\right)^{1/\rho}w(x)\, dx
\end{align*}
holds for all $(f_j)\in (L^1+L^\infty )\cap L^1_{\ell^{\rho}}(w)$
\end{thm}
In particular, we obtain the following corollary:
\begin{cor}\label{laccorvarergo}Let $(n_k)$ be a lacunary sequence, $1<r<\infty$, and $r^\prime =\frac{r}{r-1}$. Suppose that  $2\leq s<\infty$. Then, for all $f\in (L^1+L^\infty )\cap L^p_w(X)$, the weighted inequalities
$$\left\|\mathcal{V}_s^\prime f\right\|_{L^p(w)}\leq C_{p}(w)\left\|f\right\|_{L^p(w)}$$
hold if $w \in A^{\prime}_{p/r^\prime}$ and $r^\prime\leq p<\infty$, or if $w\in (A_p^{r^\prime})^\prime$ and $1<p<r$. Likewise, if  $w(x)^{r^\prime}\in A^{\prime}_1$, then the weak type inequality
\begin{align*}
\int_{\left\{x:\mathcal{V}_s^\prime f(x)>\lambda\right\}}w(x)\, dx  \leq C(w)\frac{1}{\lambda}\left\|f\right\|_{L^1(w)}
\end{align*}
holds for all $f\in (L^1+L^\infty )\cap L^1_w(X)$.
\end{cor}
\section{Applications to Ergodic $H^1$ Spaces with $A^{\prime}_p$ Weight}
\indent

In this section we apply our weighted transfer principles to obtain some new inequalities on weighted ergodic $H^1$ spaces.\\
\indent
Let now $(X,\mathscr{B} ,\mu )$ be a totally $\sigma$-finite measure space and $\{U_t:-\infty <t<\infty\}$ be a one-parameter group of measure-preserving transformations on $X$. It is assumed that $U_t(x)$ is jointly measurable in $x$ and $t$, and $U_{s+t}=U_sU_t$ for $s, t\in\mathbb{R}$. \\
\indent
We say that $w$ satisfies condition $A_p^\prime$ if there exists a constant $C$ such that for a.e. $x\in X$ 
\begin{align*}
\sup_I\left[\frac{1}{|I|}\int_Iw(U_tx)\,dt\right]\left[\frac{1}{|I|}\int_I[w(U_tx)]^{-1/(p-1)}\,dt\right]^{p-1}\leq C<\infty .
\end{align*}
The supremum is taken over all intervals $I\subset \mathbb{R}$; $|I|$ denotes the measure of $I$.\\
Here 
$w\in A_1^\prime$ if there exists a constant $C$ such that for a.e. $x$
$$\frac{1}{|I|}\int_I w(U_tx)\, dt\leq C\,{\rm{ess\, inf}}_I\,w(U_tx)$$
for all intervals $I$.
\begin{lem}\label{ergoint} Let $\mathcal{O}\subset X$ be a measurable set such that for each $x$ in $\mathcal{O}$, the set 
$$\mathcal{O}^x=\{t\in \mathbb{R}: U_tx\in \mathcal{O}\}$$
is open in $\mathbb{R}$. Then $\mathcal{O}$ can be decomposed into a disjoint union of sets, $\mathcal{O}=\cup I_i$, where the $I_i$'s are measurable sets such that, for each $x$ in $I_i$, the orbit through $x$ is a disjoint union of intervals of length between $2^k$ and $2^{k+1}$, and $k$ is an integer depending only on the set $I_i$.
\end{lem}
\begin{proof}This lemma is proved in A.~De La Torre~\cite{delatorre} when $X$ is a probability space, and one must use the fact that $X$ is totally $\sigma$-finite.
\end{proof}
\indent
Note that the set $I_i$  mentioned in Lemma~\ref{ergoint} is called an ergodic interval and the number $2^{k+1}$ is called the length of the interval.
\begin{defn}Let $\psi \in L^1({\mathbb{R}})$ and $f\in L^1(X)$, the convolution of $f$ and $\psi$ is defined as
$$f\ast \psi (x)=\int_{\mathbb{R}}f(U_{-t}x)\psi (t)\, dt.$$
\end{defn}
It is easy to observe that $f\ast \psi \in L^1(x)$ and $\|f\ast \psi \|_1\leq \|f\|_1\|\psi \|_1$. It is clear from the context in which space the norms are taken.\\
\indent
Let $\psi$ be a $C^{\infty}$ function with support in $(-1,1)$ and $L$ a positive number. For $f\in L^1(X)$ we define
$$M(L, \psi )f(x)=\sup_{|t|<\epsilon <L}|(f\ast \psi_{\epsilon})(U_tx)|$$
and
$$m(\psi )f(x)=\lim_{L\to \infty}M(L, \psi )f(x)$$
where
$$\psi_{\epsilon } (t)=\frac{1}{\epsilon}\psi (t/\epsilon ).$$
We also define a maximal operator independent of the particular $\psi$ by
$$M(L)f(x)=\sup M(L,\psi )f(x)\cdot A(\psi )^{-1}$$
where the supremum is taken over all $C^{\infty}$ functions, with support in $(-1,1)$, and $A(\psi )$ is a normalizing factor defined as
$$A(\psi )=\|\psi \|_{\infty}+\|\psi^\prime \|_{\infty}+\|\psi^{\prime\prime} \|_{\infty}.$$
We finally define 
$$Mf(x)=\lim_{L\to\infty}M(L)f(x).$$
\indent
Analogous to the maximal function characterization of real Hardy space $H^1(\mathbb{R})$ given  in C.~Fefferman and E.~Stein~\cite{fefstein}, A.~de la Torre~\cite{delatorre} defined the ergodic Hardy space $H^1(X)$ characterized by the maximal function $Mf$ as follows
$$H^1(X)=\{f\in L^1(X): Mf\in L^1(X)\}$$
with the norm
$$\|f\|_{H^1}=\|Mf\|_1.$$
\indent
Let $I$ be an ergodic interval and for $y$ fixed we consider the orbit through $y$,
$$I^y=\{t\in\mathbb{R}: U_ty\in I\}.$$
Let $I(y)$ be the interval in $I^y$ containing the origin. Then let
$$m_{I(y)}f=\frac{1}{|I(y)|}\int_{I(y)}f(U_ty)\, dt.$$
\indent
It is easy to observe that  if $z=U_ty$, $t\in I(y)$, then $m_{I(y)}=m_{I(z)}$, i.e., $m_{I(y)}$ is constant on each of the intervals forming the orbits.\\
\indent
Let now $w$ be a weight function,  we can define the weighted ergodic $H^1(X)$ spaces as follows
$$H_w^1(X)=\{f\in L_w^1(X): Mf\in L_w^1(X)\}$$
with the norm
$$\|f\|_{H_w^1}=\|Mf\|_{L_w^1}.$$
\indent
Let $f$ be a measurable function defined on $\mathbb{R}$.
$$\mathcal{V}f(x)=\left(\sum_{n=-\infty}^\infty|A_nf(x)-A_{n-1}f(x)|^s\right)^{1/s}$$
for $2\leq s<\infty$, where 
$$A_nf(x)=\frac{1}{2^n}\int_x^{x+2^n}f(y)\, dy$$
as in the previous section.\\
\indent
The following result has been obtained by the author in S. Demir~\cite{demirbsvhp}:
\begin{thm}\label{demirdyvh1l1}Let $w\in A_p$, $1\leq p<\infty$, and $s\geq 2$. Then there exists a positive constant $C$ such that
$$\|\mathcal{V}f\|_{L^1_w}\leq C\|f\|_{H^1_w}$$
for all $f\in H^1_w(\mathbb{R})$.
\end{thm}
\indent
As in the previous section consider now for a locally integrable function $f$  the two-sided variation operator 
$$\mathcal{V}^\prime f(x)=\left(\sum_{n=-\infty}^\infty|B_{2^n}f(x)-B_{2^{n-1}}f(x)|^s\right)^{1/s}$$
with
$$B_nf(x)=\frac{1}{n}\int_0^nf(U_{-t}x)\, dt.$$
\indent
Then we obtain the following result:
\begin{thm}\label{demirdyvh1l1dyn}Let $w\in A^{\prime}_p$, $1\leq p<\infty$, and $s\geq 2$. Then there exists a positive constant $C$ such that
$$\|\mathcal{V}^{\prime}f\|_{L^1_w}\leq C\|f\|_{H^1_w}$$
for all $f\in (L^1+L^\infty )\cap H^1_w(X)$.
\end{thm}
\begin{proof} This result follows by applying Theorem~\ref{wifep2} to Theorem~\ref{demirdyvh1l1} along with the maximal function characterizations of  $H^1_w(\mathbb{R})$ and  $H^1_w(X)$ spaces.
\end{proof}
\indent
Let now $f$ be a locally integrable function defined on $\mathbb{R}$, and let $(n_k)$ be a lacunary sequence. Let
$$A_{n_k}f(x)=\frac{1}{n_k}\int_0^{n_k}f(x-t)\, dt$$
and
$$\mathcal{V}_sf(x)=\left(\sum_{k=1}^\infty|A_{n_k}f(x)-A_{n_{k-1}}f(x)|^s\right)^{1/s}$$
for $2\leq s<\infty$ as in the previous section.\\
\indent
The following result has been obtained by the author in S. Demir~\cite{h1wl1w}:
\begin{thm}\label{h1wl1wvar}Let $w\in A_p$, $1\leq p<\infty$, and $s\geq 2$. Suppose that $(n_k)$ is a lacunary sequence.  Then there exists a positive constant $C$ such that
$$\|\mathcal{V}_sf\|_{L^1_w}\leq C\|f\|_{H^1_w}$$
for all $f\in H^1_w(\mathbb{R})$.
\end{thm}
Now let
$$B_nf(x)=\frac{1}{n}\int_0^nf(U_{-t}x)\, dt$$
and
$$\mathcal{V}_s^{\prime}f(x)=\left(\sum_{k=1}^\infty|B_{n_k}f(x)-B_{n_{k-1}}f(x)|^s\right)^{1/s}$$
for a locally integrable function $f$ as in the previous section. \\
\indent
When we apply Theorem~\ref{wifep2} to Theorem~\ref{h1wl1wvar} we obtain the following result:
\begin{thm}\label{h1wl1wvardyn}Let $w\in A^{\prime}_p$, $1\leq p<\infty$, and $s\geq 2$. Suppose that $(n_k)$ is a lacunary sequence.  Then there exists a positive constant $C$ such that
$$\|\mathcal{V}_s^{\prime}f\|_{L^1_w}\leq C\|f\|_{H^1_w}$$
for all $f\in (L^1+L^\infty )\cap  H^1_w(X)$.
\end{thm}

\sloppy


\begin{thebibliography}{99}
\bibitem{ha}H.~Aimer, \emph{On weighted inequalities for ergodic operators}, Studia Math. 82 (1985) 265-269. \url{http://matwbn.icm.edu.pl/ksiazki/sm/sm82/sm82115.pdf}
\bibitem{naebtag}N.~Asmar, A.~Berkon, and T. A.~Gillespie, \emph{Transference of strong type maximal inequalities by separation-preserving representations}, Amer. J. Math. 113 (1) (1991) 47-74. \url{https://www.jstor.org/stable/2374821}
\bibitem{nhasjms}N. H.~Asmar, and S. J.~Montgomery-Smith, \emph{Transference in spaces of measures}, J. Func. Analysis 165 (1999) 1-23. \url{https://www.sciencedirect.com/science/article/pii/S0022123699934055}
\bibitem{bbmcg}B.~Berkson, O.~Blasco, M. J.~Carro, and T. A.~Gillespie, 
\emph{Discretization and transference of bisublinear maximal operators}, The J. of Fourier Anal. \& Appl. 12 (4) (2006) 447-481. \url{https://link.springer.com/article/10.1007/s00041-006-6030-0}
\bibitem{apcal}A. P.~Calder\'on, 
\emph{Ergodic theory and translation-invariant operators}, Proc. Nat. Acad. Sci. USA 59 (1968) 349-353. \url{https://www.pnas.org/doi/abs/10.1073/pnas.59.2.349}
\bibitem{rrccf} R. R.~Coifman and C.~Fefferman,
\emph{Weighted norm inequalities for maximal functions and singular integrals},
Studia Math. 51 (1974) 241-250. \url{https://www1.cmc.edu/pages/faculty/MONeill/Math%20138/papers138/CoifmanFefferman.pdf}
\bibitem{rrcgw1}R. R.~Coifman, and G.~Weiss, 
\emph{Operators associated with representations of amenable groups singular integrals induced by ergodic flows, the rotation method and multipliers}, Stud. Math. 47 (1973) 285-303. \url{https://eudml.org/doc/217797}
\bibitem{rrcgw2}R. R.~Coifman, and G.~Weiss, 
\emph{Transference methods in analysis}, 
C.B.M.S. Regional Conf. Series in Math.  No. 31, AMS, Providence, Rhode Island, May 31-June 4, 1976. \url{http://www.ams.org/books/cbms/031/}
\bibitem{delatorre}A.~De La Torre, 
\emph{Ergodic $H^1$ spaces},
Bol. Soc. Mat. Mexicana 22 (1977) 10-22. \url{https://www.boletin.math.org.mx/pdf/2/22/BSMM(2).22.10-22.pdf}
\bibitem{sdem}S.~Demir, 
\emph{$H^p$ spaces and inequalities in ergodic theory},
Ph.D Thesis,  University of Illinois at Urbana-Champaign, USA, May 1999. \url{https://www.ideals.illinois.edu/items/88254}
\bibitem{demgct}S.~Demir, 
\emph{A generalization of Calder\'on transfer principle}, 
J. Comp. \& Math. Sci. 9 (5) (2018) 325-329. \url{https://www.researchgate.net/publication/326754498_A_Generalization_of_Caldern_Transfer_Principle}
\bibitem{demect}S.~Demir, 
\emph{An extension of Calder\'on Transfer Principle and its application to ergodic maximal function}, 
Asian J. of Mathematical Sciences, 4 (2) (2020) 15-18. \url{http://www.ajms.in/index.php/ajms/article/view/272/154}
\bibitem{sdemir}S.~Demir,
\emph{Inequalities for the variation operator},
Bull. of Hellenic Math Soc. 64 (2020) 92-97. \url{https://bulletin.math.uoc.gr/bulletin/vol/64/64-92-97.pdf}
\bibitem{sdemirvar}S.~Demir,
\emph{Variational inequalities for the differences of averages over lacunary sequences}, New York J. Math. 28 (2022) 1099-1111. \url{http://nyjm.albany.edu/j/2022/28-46v.pdf}
\bibitem{sdemosch1} S.~Demir, 
\emph{Oscillation inequalities on real and ergodic $H^1$ spaces},
Russian Mathematics, 67 (3) (2023) 42-52. \url{https://link.springer.com/article/10.3103/S1066369X23030039}
\bibitem{sqrfcharrnerh1} S.~Demir, 
\emph{Square function characterization of  real and ergodic $H^1$ spaces},
Russian Mathematics, 67 (4) (2023) 11-21. \url{https://link.springer.com/article/10.3103/S1066369X23040023}
\bibitem{h1wl1w}S.~Demir,
\emph{The variation operator of differences of averages over lacunary sequences maps $H_w^1(\mathbb{R})$ to$L_w^1(\mathbb{R})$}, Russian Math. 68 (5) (2024) 20-26. \url{https://link.springer.com/article/10.3103/S1066369X24700324}
\bibitem{demirbsvhp}S.~Demir,
\emph{Banach space-valued $H^p$ spaces with $A_p$ weight}, Illinois J. Math. 68 (2) (2024) 331-339. \url{https://projecteuclid.org/journals/illinois-journal-of-mathematics/volume-68/issue-2/Banach-spacevalued-Hp-spaces-with-Ap-weight/10.1215/00192082-11321393.short}
\bibitem{vlamdajump}S.~Demir,
\emph{Variation and $\lambda$-jump inequalities on $H^p$ spaces}, Russian Math. 68 (4) (2024) 12-16. \url{https://link.springer.com/article/10.3103/S1066369X24700233}
\bibitem{fefstein} C.~Fefferman and E.~Stein, 
\emph{$H^p$ spaces of several variables},
Acta Math. 129 (1972) 137-193. \url{https://projecteuclid.org/journals/acta-mathematica/volume-129/issue-none/Hp-spaces-of-several-variables/10.1007/BF02392215.full}
\bibitem{cfpw}C.~Finet, and P.~Wantiezi,
\emph{Transfer principles and ergodic theory in Orlicz spaces}, Note di mathematics, 25 (1) (2005/2006) 167-189. \url{http://siba-ese.unisalento.it/index.php/notemat/article/view/1051}
\bibitem{jgc}J.~Garcia-Cuerva,
\emph{Weighted $H^p$ spaces}, Dissertations Math. 162 (1979) 1-63.
\url{http://pldml.icm.edu.pl/pldml/element/bwmeta1.element.zamlynska-e75bca70-ae0a-4b42-b3cb-fc9e2b3cbad8}
\bibitem{jgcjrdf}J.~Garcia-Cuerva, and J. L.~Rubio de Francia,
\emph{Weighted norm inequalities and related topics}, Mathematics Studies 116,  North-Holland 1985. \url{https://www.sciencedirect.com/bookseries/north-holland-mathematics-studies/vol/116}
\bibitem{dkosz}D.~Kosz, \emph{Sharp constants in inequalities admitting the Calder\'on transference principle}, Ergodic Th. \& Dyn. Sys. 44 (6) (2024) 1597-1608.   \url{https://doi.org/10.1017/etds.2023.59}.
\bibitem{dskrlw} D. S.~Kurtz, and R. L.~Wheeden,
\emph{Results on weighted norm inequalities for multipliers},
Trans. AMS 255 (1979) 343-362. \url{https://www.ams.org/journals/tran/1979-255-00/S0002-9947-1979-0542885-8/}                                                                                                                                                                                                                                
\bibitem{bm1}B.~Muckenhoupt, \emph{Weighted norm inequalities for the Hardy maximal function}, Trans. AMS, 165 (1972) 207-226. \url{https://www.ams.org/journals/tran/1972-165-00/S0002-9947-1972-0293384-6/}
\bibitem{kp1} K.~Petersen,
\emph{Ergodic Theory}, Cambridge University Press, Cambridge, 1989. \url{https://www.cambridge.org/core/books/ergodic-theory/4F50E2830B2812125F24D4A2CE7318D0}
\end{thebibliography}
\end{document}